\documentclass[11pt, a4paper]{amsart}

\usepackage{amsmath}
\usepackage{amssymb}
\usepackage{amsthm}
\usepackage{bbm}

\usepackage{psfrag}
\usepackage{graphicx}
\usepackage{booktabs}
\usepackage{url}
\usepackage[footnotesize,nooneline,bf]{caption2}
\usepackage[mathscr]{eucal} 

\usepackage{xspace}

\graphicspath{{Pictures/}}   

\newtheoremstyle{plainNoItalics}{}{}{\normalfont}{}{\bfseries}{.}{ }{}

\newtheorem{theorem}{Theorem}[section]
\newtheorem{proposition}[theorem]{Proposition}
\newtheorem{corollary}[theorem]{Corollary}
\newtheorem{lemma}[theorem]{Lemma}
\newtheorem*{proposition*}{Proposition}
\newtheorem*{lemma*}{Lemma}
\theoremstyle{plainNoItalics}

\newtheorem*{remark}{Note}
\newtheorem*{observation}{Observation}

\newcommand{\PP}{\ensuremath{\mathcal{P}}\xspace}
\newcommand{\NP}{\ensuremath{\mathcal{NP}}\xspace}
\newcommand{\R}{\mathbbm{R}}
\newcommand{\Z}{\mathbbm{Z}}
\newcommand{\Q}{\mathbbm{Q}}
\newcommand{\card}[1]{\lvert {#1} \rvert}
\newcommand{\suchthat}{\, : \,}
\newcommand{\vones}{\boldsymbol{\mathbbm{1}}}
\newcommand{\order}[1]{\mathcal{O}\!\left({#1}\right)}
\renewcommand{\vec}[1]{{\boldsymbol {#1}}}
\newcommand{\T}{\mathrm{T}}

\newcommand{\HD}{\ensuremath{H}\xspace}
\newcommand{\faces}{\mathscr{F}}
\newcommand{\levels}{\mathscr{F}}
\newcommand{\cycles}{\mathcal{C}}
\newcommand{\PMM}{P_{M}}
\newcommand{\MMM}{\text{\textsc{MaxMM}}\xspace}
\newcommand{\inc}[1]{\vec{I}({#1})}

\DeclareMathOperator{\conv}{conv}

{\end{list}}

\newenvironment{myitemize}{%
\begin{list}{$\circ$}%
{\setlength{\topsep}{1ex}%
\setlength{\partopsep}{0mm}%
\setlength{\parskip}{1ex}%
\setlength{\parsep}{0mm}%
\setlength{\itemsep}{1ex}%
\setlength{\labelwidth}{4mm}%
\setlength{\leftmargin}{0ex}%
\addtolength{\leftmargin}{\labelwidth}%
\addtolength{\leftmargin}{\labelsep}%
\setlength{\itemindent}{0mm}}}%
{\end{list}}

\hoffset= - 0.75 in

\title{Computing Optimal Morse Matchings}
\date{8/23/2004}
\author[Joswig and Pfetsch]{Michael Joswig \and Marc E. Pfetsch}
\address{Michael Joswig, Institut f\"ur Mathematik, MA 6-2, TU Berlin,
  10623 Berlin, Germany}
\email{joswig@math.tu-berlin.de}
\address{Marc E. Pfetsch, Zuse Institute Berlin, Takustr.~7, 14195
  Berlin, Germany
}
\email{pfetsch@zib.de}
\thanks{Both authors are (partially) supported by the DFG Research Center ``Matheon''}

\begin{document}

\begin{abstract}
  Morse matchings capture the essential structural information of
  discrete Morse functions. We show that computing optimal Morse
  matchings is \NP-hard and give an integer programming formulation
  for the problem. Then we present polyhedral results for the
  corresponding polytope and report on computational results.
\end{abstract}

\maketitle

\section{Introduction}

Discrete Morse theory was developed by
Forman~\cite{For98,MR2003j:57040} as a combinatorial analog to the
classical smooth Morse theory.  Applications to questions in
combinatorial topology and related fields are numerous: e.g., Babson
et al.~\cite{MR2000a:57001}, Forman~\cite{MR2001k:57006}, Batzies and
Welker~\cite{MR2003b:13017}, and Jonsson~\cite{MR2004h:05130}.

It turns out that the topologically relevant information of a discrete
Morse function $f$ on a simplicial complex can be encoded as a
(partial) matching in its Hasse diagram (considered as a graph), the
\emph{Morse matching} of~$f$.  A matching in the Hasse diagram is
Morse if it satisfies a certain, entirely combinatorial, acyclicity
condition.  Unmatched $k$-dimensional faces are called
\emph{critical}; they correspond to the critical points of rank~$k$
of a smooth Morse function.  The total number of non-critical faces
equals twice the number of edges in the Morse matching.  The purpose
of this paper is to study algorithms which compute maximum Morse
matchings of a given finite simplicial complex.  This is equivalent to
finding a Morse matching with as few critical faces as possible.

A Morse matching $M$ can be interpreted as a discrete flow on a
simplicial complex~$\Delta$.  The flow indicates how $\Delta$ can be
deformed into a more compact description as a CW complex with one cell
for each critical face of~$M$.  Naturally one is interested in a most
compact description, which leads to the combinatorial optimization
problem described above.  This way optimal (or even sufficiently good)
Morse matchings of~$\Delta$ can help to recognize the topological type
of a space given as a finite simplicial complex.  The latter problem
is known to be undecidable even for highly structured classes of
topological spaces, such as smooth $4$-manifolds.  

Optimization of discrete Morse matchings has been studied by Lewiner,
Lopes, and Tavares~\cite{LewLT03}.  Hersh~\cite{Her03} investigated
heuristic approaches to the maximum Morse matching problem with
applications to combinatorics.  Morse matchings can also be
interpreted as pivoting strategies for homology computations;
see~\cite{Jos03}.  Furthermore, the set of all Morse matchings of a
given simplicial complex itself has the structure of a simplicial
complex; see~\cite{ChaJ01}.

Since its beginnings in Lov\'asz' proof~\cite{MR81g:05059} of the
Kneser Conjecture, combinatorial topology seeks to solve combinatorial
problems with techniques from (primarily algebraic) topology.  And,
conversely, such applications to combinatorics frequently shed light
on basic concepts in topology.  Recent years saw a still growing
influx of ideas from differential geometry to the subject.  Besides
discrete Morse theory this concerns, e.g., various notions of discrete
curvature which have been applied to problems in computational
geometry and computer graphics.  This is paralleled in the development
of a discrete differential geometry.  Here the combinatorial cores of
known phenomena in differential geometry are investigated, while
keeping a steady eye on applications in computer graphics and
mathematical physics; see, e.g., Pinkall and
Polthier~\cite{MR94j:53009}, Bobenko and Suris~\cite{MR2003d:37127}.
We firmly believe that ---via discrete Morse theory--- techniques from
combinatorial optimization are highly relevant to these topics.  Our
contribution is a first step.

The paper is structured as follows.  First we show that computing
optimal Morse matchings is \NP-hard.  This issue had been addressed
previously by Lewiner, Lopes, and Tavares~\cite{LewLT03}, but we
believe that their argument has a minor gap.  Then we give an integer
programming formulation for the problem.  The formulation consists of
two parts: one for the matching conditions and one for the acyclicity
constraints.  This turns out to be related to the acyclic subgraph
problem studied by Gr\"otschel, J\"unger, and Reinelt~\cite{GroJR85}.
We derive polyhedral results for the corresponding polytope.  In
particular, we give two different polynomial time algorithms for the
separation of the acyclicity constraints.  The paper closes with
computational results.

Like most of discrete Morse theory, also most of our results extend to
arbitrary finite CW-complexes.  We stick to the simplicial setting,
however, to simplify the presentation.

\section{Discrete Morse Functions and Morse Matchings}
\label{sec:DiscreteMorseFunctionsAndMorseMatchings}

We will first introduce discrete Morse functions as developed by
Forman~\cite{For98}. Chari~\cite{Cha00} showed that the essential
structure of discrete Morse functions is captured by so-called Morse
matchings. It turns out that this latter formulation directly leads to
a combinatorial optimization problem in which one wants to maximize
the size of a Morse matching.\smallskip

We first need some notation. Let $\Delta$ be a \emph{(finite abstract)
  simplicial complex}, i.e., a set of subsets of a finite set $V$ with
the following property: if $F \in \Delta$ and $G \subseteq F$, then $G
\in \Delta$; in other words, $\Delta$ is an independence system with
ground set~$V$. In the following we will ignore $\varnothing$ as a
member of~$\Delta$. The elements in $V$ are called \emph{vertices} and
the sets in~$\Delta$ are called \emph{faces}. The \emph{dimension} of
a face~$F$ is $\dim F := \card{F} -1$.  Let $d = \max \{ \dim F
\suchthat F \in \faces\}$ be the dimension of~$\Delta$. We often write
$i$-faces for $i$-dimensional faces.  Let~$\faces$ be the set of faces
of~$\Delta$ and let~$f_i = f_i(\Delta)$ be the number of faces of
dimension~$i \geq 0$.  The maximal faces with respect to inclusion are
called \emph{facets} and $1$-faces are called edges. The
complex~$\Delta$ is \emph{pure}, if all facets have the same
dimension. For $F$, $G \in \Delta$, we write $F \prec G$ if $F \subset
G$ and $\dim F = \dim G - 1$, i.e., ``$\prec$'' denotes the covering
relation in the Boolean lattice. The \emph{graph} of~$\Delta$ is the
(abstract) graph on~$V$ in which two vertices are connected by an edge
if there exists a $1$-face containing both vertices. Throughout this
paper we assume that $\Delta$ is \emph{connected}, i.e., its graph is
connected. This is no loss of generality since the connected
components can be treated separately.

The \emph{size} of~$\Delta$ is defined as the coding length of its
face lattice, i.e., if $\Delta$ has $n$ faces, then
$\operatorname{size}\Delta = \order{n \cdot d \cdot \log n}$.
Statements about the complexity of algorithms in the subsequent
sections are always with respect to this notion of size.

A function $f : \Delta \rightarrow \R$ is a \emph{discrete Morse
  function} if for every $G \in \Delta$ the sets
\begin{equation}
\label{eq:MorseCond}
\{ F : F \prec G,\; f(G) \leq f(F) \}\quad \text{and}\quad
\{ H : G \prec H,\; f(H) \leq f(G) \}
\end{equation}
both have cardinality at most~$1$. The first set includes the faces
covered by face~$G$ which are not assigned a lower value than~$G$,
while the second set includes the faces covering~$G$ which are not
assigned a higher value. The face~$G$ is \emph{critical} if both sets
have cardinality~$0$. A simple example of a discrete Morse function
can be obtained by setting $f(F) = \dim F$ for every $F \in \Delta$.
With respect to this function every face is critical.

Discrete Morse functions are interesting because they can be used to
deform a simplicial complex into a (smaller) CW-complex that has a
cell for each critical face; see
Section~\ref{sec:PropertiesOfMorseMatchings}.\smallskip

Consider the \emph{Hasse diagram}~$\HD = (\faces, A)$ of $\Delta$,
that is, a directed graph on the faces of~$\Delta$ with an arc $(G,F)
\in A$ if $F \prec G$; note that the arcs lead from higher to lower
dimensional faces. Let $M \subset A$ be a matching in~\HD, i.e., each
face is incident to at most one arc in~$M$. Let~$\HD(M)$ be the
directed graph obtained from~\HD by reversing the direction of the
arcs in~$M$. Then~$M$ is a \emph{Morse matching} of~$\Delta$
if~$\HD(M)$ does not contain directed cycles, i.e., is acyclic (in the
directed sense).

Chari~\cite{Cha00} observed the following relation between Morse
functions and Morse matchings. Let~$f$ be a discrete Morse function
and let $M$ be the set of arcs $(G,F) \in A$ such that $f(G) \leq
f(F)$, i.e., $f$ is not decreasing on these arcs. A simple proof shows
that at most one of the sets in~\eqref{eq:MorseCond} can have
cardinality one; see Chari~\cite{Cha00}. This shows that~$M$ is a
matching. Since the order given by~$f$ can be refined to a linear
ordering of the faces of~$\Delta$, the directed graph~$\HD(M)$ is in
fact acyclic and therefore a Morse matching. To construct a discrete
Morse function from a Morse matching, compute a linear ordering
extending~$\HD(M)$ (which is acyclic) and then number the faces
consecutively in the reverse order.

Although we loose the concrete numbers attached to the faces when
going from a discrete Morse function~$f$ to the corresponding Morse
matching~$M$, we do not loose the information about critical faces:
Critical faces of $f$ are exactly the unmatched faces of~$M$.  Hence,
by maximizing~$\card{M}$ we minimize the number of critical faces
of~$f$. In fact, the number of critical faces is $\card{\faces} - 2\,
\card{M}$. For $0 \leq j \leq d$, let $c_j = c_j(M)$ be the number of
critical faces of dimension~$j$ and let $c(M)$ be the total number of
critical faces.

It seems helpful to briefly describe the case of Morse matchings for a
one-dimensional simplicial complex~$\Delta$. Then $\Delta$ represents
the incidences of a graph~$G$. A Morse matching~$M$ of~$\Delta$
matches edges with nodes of~$G$. Let $\tilde{G}$ be the following
oriented subgraph of~$G$: take all edges which are matched in~$M$ and
orient them towards its matched node. Since~$M$ is a matching, this
construction is well defined and the in-degree of each node is one.
The acyclicity property shows that $\tilde{G}$ contains no directed
cycles and hence is a branching, i.e., the underlying graph is a
forest and each (weakly) connected component has a unique root.
Therefore, the Morse matchings on a graph $G$ are in one-to-one
correspondence with orientations of subgraphs of~$G$ which are
branchings.

Generalizing this idea, Lewiner, Lopes, and Tavares~\cite{LewLT03}
developed a heuristic for computing optimal Morse matchings. This
heuristic computes maximum Morse matchings for combinatorial
$2$-manifolds, i.e., Morse matchings with maximal cardinality.
However, for general simplicial complexes this problem is \NP-hard,
see Section~\ref{sec:Hardness}.

\section{Properties of Morse Matchings}
\label{sec:PropertiesOfMorseMatchings}

In this section we briefly review some important properties of Morse
matching which we need in the sequel.

Let $F$ be a facet of $\Delta$ and let $G$ be a facet of $F$, which is
not contained in any other facet of $\Delta$. The operation of
transforming $\Delta$ to $\Delta \setminus \{F, G\}$ is called a
\emph{simplicial} or \emph{elementary collapse}. We will simply use
collapse in the following.

\begin{proposition}[Forman~\cite{For98}]\label{prop:subcomplex_collapse}
  Let $\Delta$ be a simplicial complex and $\Sigma$ a subcomplex of
  $\Delta$. Then there exists a sequence of collapses from $\Delta$ to
  $\Sigma$ if and only if there exists a discrete Morse function such
  that $\Delta \setminus \Sigma$ contains no critical face.
\end{proposition}

Forman~\cite{For98} also proved the following result, which describes
one of the most interesting features of Morse matchings:

\begin{theorem}\label{thm:CWComplex}
  Let $\Delta$ be a simplicial complex and $M$ be a Morse matching on
  $\Delta$. Then $\Delta$ is homotopy equivalent to a CW-complex
  containing a cell of dimension~$i$ for each critical face of
  dimension~$i$.
\end{theorem}

We refer to Munkres~\cite{Mun84} for more information on CW-complexes.
By this result we can hope for a compact representation of the
topology of~$\Delta$ (up to homotopy) by computing a Morse matching
with few critical faces. This is the main motivation for the
combinatorial optimization problem studied in this paper.

Let~$K$ be a field and let $\beta_j = \beta_j(K)$ to be the Betti
number for dimension~$j$ over~$K$ for~$\Delta$; see again
Munkres~\cite{Mun84} for details. Forman~\cite{For98} proved the
following bounds on the number of critical faces $c_j$ of a Morse
matching~$M$:

\begin{theorem}[\emph{Weak Morse inequalities}]\label{thm:weak_Morse_inequalities}
  Let $K$ be a field, $\Delta$ be a simplicial complex, and $M$ a
  Morse matching for $\Delta$. We have
  \begin{equation}\label{eq:morse_inequalities}
  c_j \geq \beta_j \qquad \text{for all }j = 0, \dots, d
  \end{equation}
  and
  \begin{equation}\label{eq:morse_equation}
  c_0 - c_1 + c_2 - \dots + (-1)^d c_d = \beta_0 - \beta_1 + \beta_2
  - \dots + (-1)^d \beta_d.
  \end{equation}
\end{theorem}

The Betti numbers over~$\Q$ and finite fields can easily be obtained
in polynomial time (in the size of~$\Delta$), by computing the ranks
of the boundary matrices for each dimension. Although harder to
compute (see Iliopoulos~\cite{Ili89}), the homology over~$\Z$ can be
used to choose among the finite fields or~$\Q$, in order to obtain the
strongest form of the Morse
inequalities~\eqref{thm:weak_Morse_inequalities}.

\section{Hardness of Optimal Morse Matchings}
\label{sec:Hardness}

In this section we prove \NP-hardness of the problem to compute a
maximum Morse matching, i.e., to find a Morse matching~$M$ with
maximal cardinality. As we saw previously, this is equivalent to
minimize the number of critical faces.

We want to reduce the following \emph{collapsibility problem},
introduced by E\v{g}ecio\v{g}lu and Gonzalez~\cite{EgeG96}, to the
problem of finding an optimal Morse matching: Given a connected pure
$2$-dimensional simplicial complex $\Delta$, which is embeddable in
$\R^3$ and an integer~$k$, decide whether there exists a subset
$\mathcal{K}$ of the facets of $\Delta$ with $\card{\mathcal{K}} \leq
k$ such that there exists a sequence of collapses which transforms
$\Delta \setminus \mathcal{K}$ to a $1$-dimensional complex.
E\v{g}ecio\v{g}lu and Gonzalez proved that this collapsibility problem
is strongly \NP-complete. Using
Proposition~\ref{prop:subcomplex_collapse}, this results reads in
terms of discrete Morse theory:

\begin{theorem}\label{thm:erasable}
  Given a connected pure $2$-dimensional simplicial complex~$\Delta$,
  which is embeddable in~$\R^3$, and a nonnegative integer~$k$, it is
  \NP-complete in the strong sense to decide whether there exists a
  Morse matching with at most $k$ critical $2$-faces.
\end{theorem}

When~$k$ is fixed, we can try all possible sets $\mathcal{K}$ of size
at most~$k$ and then decide whether the resulting complex is
collapsible to a $1$-dimensional complex in polynomial time. Therefore
we let $k$ be part of the input.

We need the following construction. Consider a Morse matching~$M$ for
a simplicial complex~$\Delta$, with $\dim \Delta \geq 1$. Let
$\Gamma(M)$ be the graph obtained from the graph of~$\Delta$ by
removing all edges ($1$-faces) matched with $2$-faces.  Note that
$\Gamma(M)$ contains all vertices of~$\Delta$.

\begin{lemma}\label{lemma:connected}
  The graph $\Gamma(M)$ is connected.
\end{lemma}

\begin{proof}
  Without loss of generality we assume that $\dim \Delta \geq 2$.
  Otherwise $\Gamma(M)$ coincides with the graph of~$\Delta$, which
  is connected (recall that $\Delta$ is connected).
  
  Suppose that $\Gamma(M)$ is disconnected. Let $N$ be the set of
  nodes in a connected component of $\Gamma(M)$, and let $C$ be the
  set of \emph{cut edges}, that is, edges of~$\Delta$ with one
  vertex in $N$ and one vertex in its complement. Since~$\Delta$ is
  connected, $C$ is not empty. By definition of $\Gamma(M)$, each edge
  in~$C$ is matched to a unique $2$-face.
  
  Consider the directed subgraph~$D$ of the Hasse diagram consisting
  of the edges in~$C$ and their matching $2$-faces. The standard
  direction of arcs in the Hasse diagram (from the higher to the
  lower dimensional faces) is reversed for each matching pair
  of~$M$, i.e., $D$ is a subgraph of $\HD(M)$.
  
  We construct a directed path in~$D$ as follows. Start with any node
  of~$D$ corresponding to a cut edge~$e_1$. Go to the node of~$D$
  determined by the unique $2$-face~$\tau_1$ to which~$e_1$ is matched
  to. Then~$\tau_1$ contains at least one other cut edge~$e_2$,
  otherwise $e_1$ cannot be a cut edge. Now iteratively go to~$e_2$,
  then to its unique matching $2$-face~$\tau_2$, choose another cut
  edge~$e_3$, and so on. We observe that we obtain a directed path
  $e_1, \tau_1, e_2, \tau_2, \dots$ in~$D$, i.e., the arcs are
  directed in the correct direction.
  
  Since we have a finite graph at some point the path must arrive at a
  node of~$D$ which we have visited already. Hence, $D$ (and therefore
  also~$\HD(M)$) contains a directed cycle, which is a contradiction
  since $M$ is a Morse matching.
\end{proof}

Now pick an arbitrary node~$r$ and any spanning tree of $\Gamma(M)$
and direct all edges away from~$r$. This yields a maximum Morse
matching on~$\Gamma(M)$; see end of
Section~\ref{sec:DiscreteMorseFunctionsAndMorseMatchings}. It is easy
to see that replacing the part of~$M$ on $\Gamma(M)$ with this
matching yields a Morse matching. This Morse matching has only one
critical vertex (the root~$r$). Note that every Morse matching
contains at least one critical vertex; this can be seen from the Morse
inequalities Theorem~\ref{thm:weak_Morse_inequalities} or directly by
observing that we get a directed cycle in the Hasse diagram if every
vertex is matched to an edge. Furthermore, the total number of
critical faces can only decrease, since we computed an optimal Morse
matching on $\Gamma(M)$. The number of critical $i$-faces for $i \geq
2$ stays the same. We have thus proved the following, which is also
implicit in Forman~\cite{For98}.

\begin{corollary}\label{cor:critical_vertex}
  Let~$M$ be a Morse matching on $\Delta$. Then we can compute a Morse
  matching~$M'$ in polynomial time which has exactly one critical
  vertex and the same number of critical faces of dimension $2$ or
  higher as~$M$, such that $c(M') \leq c(M)$.
\end{corollary}

We can now prove the hardness result.

\begin{theorem}\label{thm:np-complete}
  Given a simplicial complex~$\Delta$ and a nonnegative integer~$c$,
  it is strongly \NP-complete to decide whether there exists a Morse
  matching with at most $c$ critical faces, even if~$\Delta$ is
  connected, pure, $2$-dimensional, and can be embedded in~$\R^3$.
\end{theorem}

\begin{proof}
  Clearly this problem is in \NP. So let~$(\Delta,\,k)$ be an input
  for the collapsibility problem. We claim that there exists a Morse
  matching with at most~$k$ critical $2$-faces if and only if there
  exists a Morse matching with at most~$g(k) := 2 (k + 1) -
  \chi(\Delta)$ critical faces altogether. Here, $\chi(\Delta) =
  \beta_0 - \beta_1 + \dots + (-1)^d \beta_d$ is the Euler
  characteristic, which can be computed in polynomial time; see
  Section~\ref{sec:PropertiesOfMorseMatchings}. Hence $g$ is a
  polynomial-time computable function. Using
  Theorem~\ref{thm:erasable} then finishes the proof.
  
  So assume that~$M$ is a Morse matching on~$\Delta$ with at most~$k$
  critical $2$-faces. We use Corollary~\ref{cor:critical_vertex} to
  compute a Morse matching~$M'$, in polynomial time, such that
  $c_0(M') = 1$, $c_2(M') = c_2(M)$, and $c(M') \leq c(M)$. By the
  Morse equation of Theorem~\ref{thm:weak_Morse_inequalities}, we have
  $c_1(M') = c_2(M') + 1 - \chi(\Delta)$. Since $c(M') = c_0(M') +
  c_1(M') + c_2(M')$ it follows that
  \begin{equation}\label{eq:critical_triangles}
  c_2(M) = c_2(M') = \tfrac{1}{2}(c(M') + \chi(\Delta))) - 1.
  \end{equation}
  Solving for $c(M')$, it follows that $M'$ has at most $2 (k + 1) -
  \chi(\Delta)$ critical faces altogether.
  
  Conversely, assume that there exists a Morse matching~$M$ with at
  most $g(k)$ critical faces. Computing $M'$ as above, we obtain
  by~\eqref{eq:critical_triangles}, that
  \[
  c_2(M) = c_2(M') \leq \tfrac{1}{2}(g(k) + \chi(\Delta))) - 1 = k,
  \]
  which completes the proof.
\end{proof}

Lewiner, Lopes, and Tavares~\cite{LewLT03} showed that it is \NP-hard
to compute an optimal Morse matching with exactly one critical vertex.
They claimed it for the general case, but we do not see an argument
similar to Lemma~\ref{lemma:connected} above.

Since there exists a Morse matching with at most $c$ critical faces if
and only if there exists a Morse matching of size at least
$\tfrac{1}{2} (\card{\faces} - c)$, we proved:

\begin{corollary}
  Let $\Delta$ be as in Theorem~\ref{thm:np-complete} and $m$ be a
  nonnegative integer. Then it is \NP-complete in the strong sense to
  decide whether there exists a Morse matching of size at least~$m$.
\end{corollary}

We do not know about the complexity status for this problem with $m$
fixed.

E\v{g}ecio\v{g}lu and Gonzalez~\cite{EgeG96} even proved that the
collapsibility problem is as hard to approximate as the set covering
problem. In particular, the collapsibility problem cannot be
approximated better than within a logarithmic factor in polynomial
time, unless $\PP = \NP$. Using this, Lewiner, Lopes, and
Tavares~\cite{LewLT03} claimed that the problem to compute a Morse
matching minimizing the number of critical faces is hard to
approximate. However, the function~$g$ used in the proof above is not
``approximation preserving'' and we do not see how the
non-approximability result carries over.

Similarly, the problem to approximate the size of a Morse matching
seems to be open.

\section{An IP-Formulation}
\label{sec:ip_formulation}

In this section we introduce an integer programming formulation for
the problem to compute a Morse matching of maximal size. It can easily
be extended to arbitrary weights.

We use the following notation. We depict vectors in bold font. Let
$\vec{e}_i$ be the $i$th unit vector and let~$\vones$ be the vector of
all ones. For any vector $\vec{x} \in \R^n$ and $I \subseteq \{1,
\dots, n\}$ we define
\[
\vec{x}(I) := \sum_{i \in I} x_i.
\]
Furthermore, for $S \subseteq \{1, \dots n\}$, $\inc{S} \in \R^n$
denotes the incidence vector of~$S$.

For a node $v$ in a directed graph, let $\delta(v)$ be the arcs
incident to $v$, i.e., the arcs having $v$ as one of their endnodes.
For a subset $A' \subseteq A$, we denote by $N(A')$ the nodes incident
to at least one arc in $A'$.  Throughout this article, all directed or
undirected cycles are assumed to be \emph{simple}, i.e., without node
repetitions.

For ease of notation, we consider the Hasse diagram~\HD as directed or
undirected depending on the context; we will explicitly use
\emph{directed} when we refer to the directed version.

We split $\HD$ into $d$ levels $\HD_0 = (\levels^0, A_0)$, \dots,
$\HD_{d-1} = (\levels^{d-1}, A_{d-1})$, where~$\HD_i$ denotes the
level of the Hasse diagram between faces of dimension~$i$ and~$i+1$.
Then $A$ is the disjoint union of $A_0, \dots, A_{d-1}$ and
$\levels^{i-1} \cap \levels^i$ consists of the faces of dimension~$i$.
Recall that the arcs in the Hasse diagram are directed from the higher
to the lower dimensional faces.\smallskip

Let $M \subset A$ be a Morse matching of~$\Delta$. By definition, its
incidence vector~$\vec{x} = \inc{M} \in \{0,1\}^A$ satisfies the
\emph{matching inequalities}:
\begin{equation}
  \label{eq:matching_inequalities}
  \vec{x}(\delta(F)) \leq 1 \quad  \forall\; F \in \faces.
\end{equation}

\begin{figure}
  \centering
  \includegraphics[height=2.5cm]{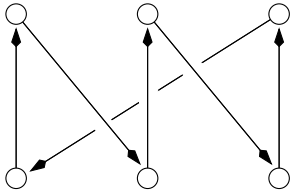}
  \caption{Example for a directed cycle of size $6$; at least three
    arcs with reversed orientation (pointing ``up'') are necessary to
    close a $6$-cycle in the Hasse diagram of a simplicial complex.}
  \label{fig:6cycle}
\end{figure}

Now assume that for some $M \subseteq A$ there exists a directed
cycle~$D$ in $\HD(M)$. Then in $D$ ``up'' and ``down'' arcs alternate;
see for example Figure~\ref{fig:6cycle}. In particular, the size of
$D$ is always even. Hence, $\tfrac{1}{2}\card{D}$ arcs are contained
in~$M$, i.e., are reversed in~$\HD(M)$. Furthermore we note:

\begin{observation}
  Let~$M \subset A$ be a matching. If $D$ is a directed cycle in
  $\HD(M)$, the edges in~$D$ can only belong to one level~$H_i$ ($i
  \in \{0, \dots, d-1\}$), i.e., $\{\dim F \suchthat F \in N(D)\} =
  \{i, i+1\}$.
\end{observation}

\noindent Putting these arguments together we obtain: If $M$ is acyclic,
$\vec{x} = \inc{M}$ satisfies the following \emph{cycle inequalities}:
\begin{equation}
  \label{eq:cycle_inequalities}
  \vec{x}(C) \leq \tfrac{1}{2} \card{C} - 1 \quad
  \forall\; C \in \cycles_i,\; i=1,\dots, d-1,
\end{equation}
where $\cycles_i$ are the cycles in $H_i$.

Conversely, it is easy to see that every $\vec{x} \in \{0,1\}^A$ which
fulfills inequalities~\eqref{eq:matching_inequalities}
and~\eqref{eq:cycle_inequalities} is the incidence vector of a Morse
matching. Hence, we arrive at the following IP formulation for the
problem to find a maximum Morse matching:
\[
  \begin{array}{l@{\quad}l@{}r@{\;}ll}
  (\MMM)& \max        & \vones^\T \vec{x}\\[1ex]
  ~     & \text{s.t.} & \vec{x}(\delta(F)) & \leq 1 & \forall\; F \in \faces\\[0.5ex]
  ~     &             & \vec{x}(C)         & \leq \tfrac{1}{2} \card{C} - 1
  & \forall\; C \in \cycles_i,\; i = 1, \dots, d-1\\[0.5ex]
  ~     &             & \vec{x}            & \in \{0,1\}^A.
  \end{array}
\]
We define the
corresponding polytope:
\[
\PMM = \conv \big\{ \vec{x} \in \{0,1\}^A \suchthat \vec{x} \text{
  satisfies}~\eqref{eq:matching_inequalities} \text{
  and}~\eqref{eq:cycle_inequalities} \big\}.
\]

Let $M$ be a Morse matching and~$\vec{x} = \inc{M}$ be its incidence
vector. Then $F \in \faces$ is a critical face with respect to~$M$ if
and only if it is unmatched by~$M$, i.e., $\vec{x}(\delta(F)) = 0$.
Hence, the total number of critical faces is:
\begin{equation}\label{eq:critical_face_equation}
c(M) = \sum_{F \in \faces} \bigl( 1 - \sum_{a \in \delta(F)} x_a \bigr) = \card{\faces} - 2
\, \sum_{a \in A} x_a = \card{\faces} - 2\, \vones^\T \vec{x},
\end{equation}
since every arc is incident to exactly two nodes. Using this formula
one can easily switch between the number of critical faces and the
number of arcs in a Morse matching.

The LP relaxation of \MMM can be strengthened by using the weak Morse
inequalities~\eqref{eq:morse_inequalities} of
Theorem~\ref{thm:weak_Morse_inequalities}.
Applying~\eqref{eq:critical_face_equation}, this yields the following
\emph{Betti inequality} for dimension~$i$:
\begin{equation}\label{eq:betti_inequality}
\sum_{F: \dim F = i} \bigl( 1 - \sum_{a \in \delta(F)} x_a
\bigr) \geq \beta_i
\qquad \Leftrightarrow \qquad
\sum_{F: \dim F = i} \; \sum_{a \in \delta(F)} x_a \leq f_i - \beta_i.
\end{equation}
Observe that we can choose the field in
Theorem~\ref{thm:weak_Morse_inequalities} to employ the Morse
inequalities in their strongest form.

\begin{remark}
  The cycle inequalities~\eqref{eq:cycle_inequalities} are similar to
  the cycle inequalities for the acyclic subgraph problem (ASP); see
  J\"unger~\cite{Jun85}, and Gr\"otschel, J\"unger, and
  Reinelt~\cite{GroJR85}. The separation problem
  for~\eqref{eq:cycle_inequalities}, however, is more complicated than
  the corresponding problem for ASP; see
  Section~\ref{sec:separation_cycle_inequalities}.
  
  Furthermore, there is a similarity to the relation between the ASP
  and the linear ordering problem (see Reinelt~\cite{Rei85}, and
  Gr\"otschel, J\"unger, and Reinelt~\cite{GroJR84}): an alternative
  formulation for our problem can be obtained by adding matching
  inequalities to the linear ordering formulation; this directly
  models discrete Morse functions as linear orderings of the faces.
  Since this formulation is based on the relation between faces, it
  leads to quadratically many variables in the number of faces;
  therefore we have opted for the above formulation, at the cost of
  having to solve the separation problem for the cycle inequalities;
  see Section~\ref{sec:separation_cycle_inequalities}.
\end{remark}

\subsection{Facial Structure of $\boldsymbol{\PMM}$}

\begin{figure}
  \mbox{}\hspace{10ex}
  \includegraphics[height=2.5cm]{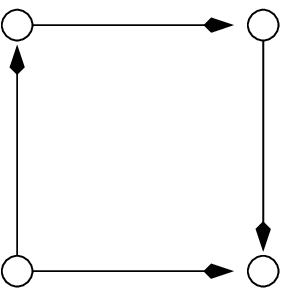}
  \hspace{15ex}
  \includegraphics[width=2.5cm]{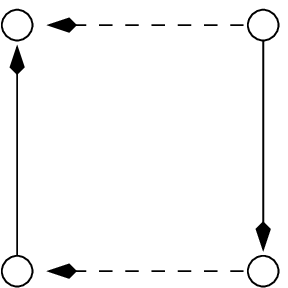}
  \mbox{}\hspace{10ex}
  \caption{Example for a non-monotone behavior of acyclic
    matchings. The directed graph on the right, obtained from the left
    graph by reversing the dashed arcs, is acyclic. However, if the
    top arc is set to its original orientation, the graph is not
    acyclic anymore. This shows that subsets of acyclic matchings are
    not necessarily acyclic.}
  \label{fig:monotone}
\end{figure}

It is easy to see that $\PMM$ is a full dimensional polytope and $x_a
\geq 0$ defines a facet for every $a \in A$. Furthermore, $\PMM$ is
monotone, since every subset of a Morse matching is a Morse matching.
It is well known that this implies that every facet defining
inequality $\vec{\alpha}^\T \vec{x} \leq \beta$ not equivalent to the
non-negativity inequalities fulfills: $\vec{\alpha} \geq 0$, $\beta >
0$; see Hammer, Johnson, and Peled~\cite{HamJP75}.

Interestingly, if we generalize Morse matchings to acyclic matchings
for arbitrary graphs, the collection of such acyclic matchings is not
necessarily monotone anymore; see the example in
Figure~\ref{fig:monotone}.  Therefore, the structure of the
generalized problem is likely to be more complicated.

We have the following two results:

\begin{proposition}
  The matching inequalities $\vec{x}(\delta(F)) \leq 1$ define facets
  of $\PMM$ for $F \in \mathcal{F}$, except if $\card{\delta(F)} = 1$,
  i.e., $F$ is a vertex.
\end{proposition}
\begin{proof}
  Let $F$ be a face with $\card{\delta(F)} > 1$ (note that
  $\card{\delta(F)} = 0$ does not occur). We can assume that $A =
  \{a_1, \dots, a_k, a_{k+1}, \dots, a_m\}$, where $\delta(F) = \{a_1,
  \dots, a_k\}$. For $i = k+1, \dots, m$, observe that~$a_i$ cannot be
  adjacent to every arc in $\delta(F)$: since $\card{\delta(F)} > 1$,
  $a_i$ would either be incident to at least two nodes of the same
  dimension or to two nodes whose dimensions are two apart, which is
  impossible. Therefore, choose $p(i) \in \{1, \dots, k\}$ such that
  $a_i$ and $a_{p(i)}$ are not adjacent. It follows that $\vec{e}_i +
  \vec{e}_{p(i)} \in \PMM$.  Then
  \[
  \vec{e}_1, \dots, \vec{e}_k, \vec{e}_{k+1} + \vec{e}_{p(k+1)},
  \dots, \vec{e}_{m} + \vec{e}_{p(m)}
  \]
  are affinely independent and fulfill $\vec{x}(\delta(F)) = 1$.
\end{proof}
It follows that the inequalities~$x_a \leq 1$, $a \in A$, never define
facets.
\begin{theorem}\label{prop:cycle_facets}
  The cycle inequalities~\eqref{eq:cycle_inequalities} define facets of
  $\PMM$.
\end{theorem}
\begin{proof}
  We extend the corresponding proof by J\"unger~\cite{Jun85} for the
  ASP.\smallskip
  
  Let $C$ be a cycle in~$\HD$. Without loss of generality
  assume that $A = \{a_1, \dots, a_k, a_{k+1}, \dots, a_m\}$, where $C
  = (a_1, \dots, a_k)$ and $k$ is even. We will construct affinely
  independent feasible vectors $\vec{v}_1, \dots, \vec{v}_k,
  \vec{v}_{k+1}, \dots, \vec{v}_m$ satisfying the cycle inequality
  corresponding to~$C$ with equality.
  
  Let $C_1 = \{a_1, a_3, \dots, a_{k-1}\}$ and $C_2 = \{a_2, a_4,
  \dots, a_k\}$.  Hence $C_1$ and $C_2$ are the ``up'' and ``down''
  arcs in $C$. 

  Define
  \[
  \vec{v}_i = 
  \begin{cases}  
    \inc{C_1 \setminus \{a_i\}} & \text{if }a_i \in C_1\\
    \inc{C_2 \setminus \{a_i\}} & \text{if }a_i \in C_2
  \end{cases}
  \qquad
  \text{for }i = 1, \dots, k.
  \]
  Hence, for $i = 1, \dots, k$ we have $\vec{v}_i(C) = \tfrac{k}{2}-1$.

  For $i = k+1, \dots, m$, consider $a_i = \{u,v\} \notin C$. We have
  four cases.

  \begin{figure}
    \centering
    \psfrag{P1}{$P_1$}
    \psfrag{P2}{$P_2$}
    \psfrag{u}{$u$}
    \psfrag{v}{$v$}
    \psfrag{C1}{$\tilde{C}_1$}
    \psfrag{C2}{$\tilde{C}_2$}
    \includegraphics[height=2.75cm]{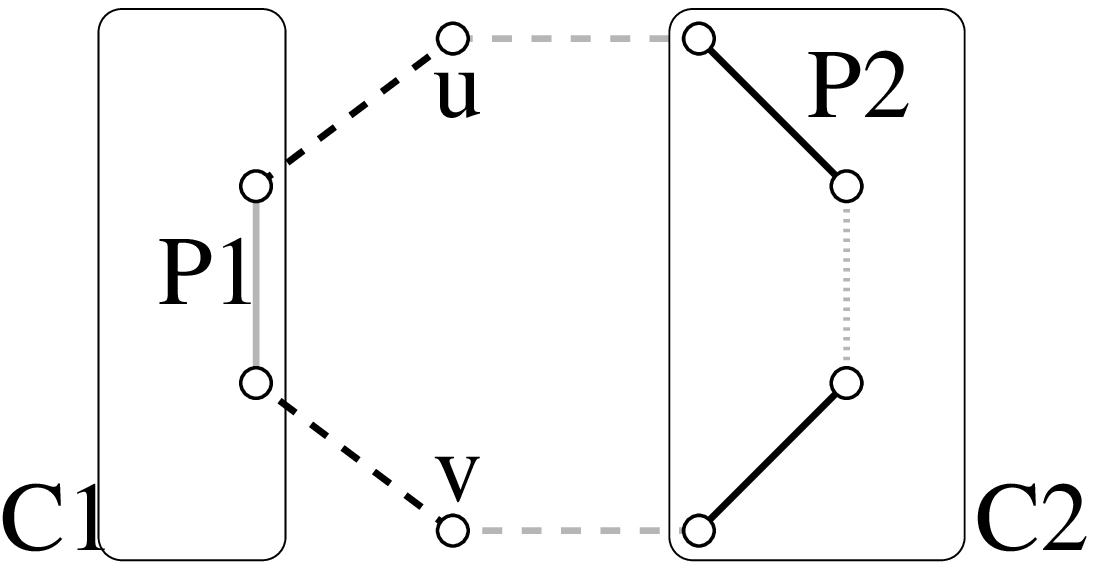}
    \hspace{10ex}
    \includegraphics[height=2.75cm]{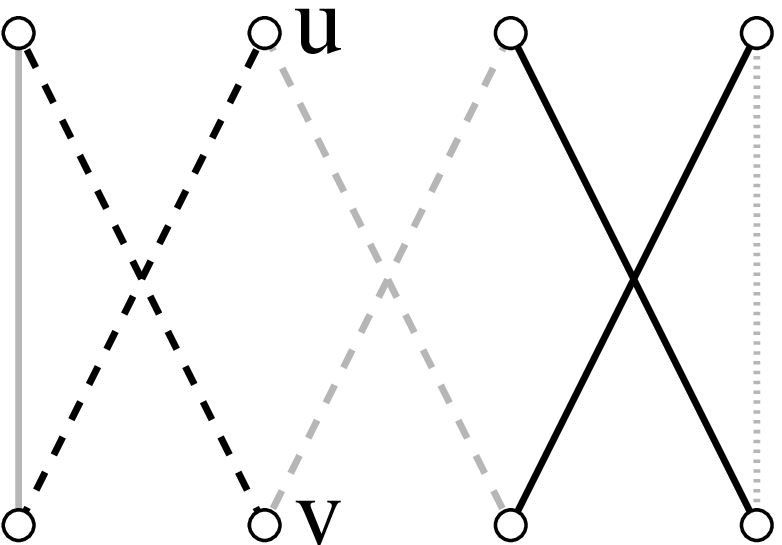}
    \caption{Illustration of the first case in the proof of
      Theorem~\ref{prop:cycle_facets}. The sets $P_1$ and $P_2$ are
      shown by continuous lines. The edges in $C_1$ are drawn gray and
      hence $P_1 \subset C_1$; edges in $C_2$ are drawn black. The
      dashed edges incident to $u$ and $v$ are not considered. The
      right hand side shows the graph embedded in the Hasse diagram.}
    \label{fig:cyclefacet}
  \end{figure}
  
  \begin{myitemize}
  \item[$\triangleright\; u,v \in N(C)$:] Let $\tilde{C} := C
    \setminus \big(\delta(u) \cup \delta(v)\big)$. We have that
    $\card{\tilde{C}} = k - 4$ (since there exist no odd cycles) and
    $\tilde{C}$ splits into two odd nonempty parts $\tilde{C}_1$ and
    $\tilde{C}_2$, which are both paths. Let $k_1 :=
    \card{\tilde{C}_1}$ and $k_2 := \card{\tilde{C}_2}$; $k_1$ and
    $k_2$ are odd, since $u$ and $v$ are on opposite sides of the
    bipartition. We choose a subset $P_1 \subset \tilde{C}_1$ by
    taking every second arc in order to get $\card{P_1} = \tfrac{k_1 +
      1}{2}$; similarly we choose $P_2 \subset \tilde{C}_2$ with
    $\card{P_2} = \tfrac{k_2 + 1}{2}$. By construction either $P_i
    \subset C_1$ or $P_i \subset C_2$ and either $P_i \cap C_2 =
    \varnothing$ or $P_i \cap C_1 = \varnothing$ for $i = 1,2$.  An
    easy calculation shows that $\card{P_1 \cup P_2} =
    \tfrac{k}{2}-1$. See Figure~\ref{fig:cyclefacet} for an
    illustration of this case. Then define $\vec{v}_i := \inc{P_1 \cup
      P_2 \cup \{a_i\}}$.
  \item[$\triangleright\; u \notin C,\; v \in C$:] Here we define
    $\vec{v}_i := \inc{C_1 \setminus \delta(v) \cup \{a_i\}}$.
  \item[$\triangleright\; u \in C,\; v \notin C$:] Define
    $\vec{v}_i := \inc{C_1 \setminus \delta(u) \cup \{a_i\}}$.
  \item[$\triangleright\; u, v \notin C$:] Choose any $a \in C_1$ and define
    $\vec{v}_i := \inc{C_1 \setminus \{a\} \cup \{a_i\}}$.
  \end{myitemize}
  It is easy to check in each case that $\vec{v}_i \in \PMM$ and that
  $\vec{v}_i(C) = \tfrac{k}{2}-1$.
  
  It can be shown that the $m$ vectors $\vec{v}_1, \dots, \vec{v}_m$
  are affinely independent, which concludes the proof.
\end{proof}

The separation problem for the cycle inequalities is discussed in the
next section.

\subsection{Separating the Cycle Inequalities}
\label{sec:separation_cycle_inequalities}

Of course, there are exponentially many cycle
inequalities~\eqref{eq:cycle_inequalities}. Hence we have to deal with
the separation problem for these inequalities.

For the separation problem, we can assume that we are given $\vec{x}^*
\in [0,1]^A$, which satisfies all matching
inequalities~\eqref{eq:matching_inequalities}. We consider the
separation for each graph~$H_i$ in turn, $i = 0, \dots, d-1$. The
problem is to find an  undirected cycle~$C$ in~$H_i$ such that
\[
\vec{x}^*(C) > \tfrac{1}{2} \card{C} - 1
\]
or conclude that no such cycle exists. In the next sections we
describe two methods to solve this separation problem in polynomial
time.

\subsubsection{Undirected Shortest Path with Conservative Weights}
\label{sec:undirected_shortest_paths}

A usual trick to solve the above separation problem is to apply an
affine transformation and obtain a shortest cycle problem. The
transformation suitable for our needs is $\vec{x}' = \tfrac{1}{2}
\vones - \vec{x}$, which yields:
\[
\vec{x}(C) \leq \tfrac{1}{2} \card{C} - 1
\qquad \Leftrightarrow \qquad
\vec{x}'(C) \geq 1.
\]
The separation problem can now be solved as follows: compute a
shortest cycle in~$H_i$ with respect to the weights $\tfrac{1}{2}
\vones - \vec{x}^*$. If its weight is at most~$1$, this cycle yields a
violated cycle inequality, otherwise no such cycle exists.

However, the weights can be negative and we have to rule out negative
cycles in order to apply polynomial time methods from the literature;
that is, we want the weights to be \emph{conservative}.

\begin{lemma}
  There exists no cycle of negative weight in~$H_i$ with respect
  to~$\tfrac{1}{2} \vones - \vec{x}^*$, for $0 \leq i \leq d-1$.
\end{lemma}
\begin{proof}
  Let $C = (a_1, \dots, a_k)$ be a cycle in~$H_i$ and let $F_1, \dots,
  F_k$ be the faces that are visited by~$C$. Recall the $\vec{x}^*$
  satisfies the matching inequalities. We obtain
  \begin{equation}\label{eq:cycle_double_count}
  \sum_{j = 1}^k \sum_{a \in \delta(F_i) \cap C} x_a^* \; = \; 2
  \sum_{a \in C} x_a^* \; = \; 2\, \vec{x}^*(C),
  \end{equation}
  since each edge weight is counted twice in the first term. Applying
  the Matching inequalities~\eqref{eq:matching_inequalities} on the
  left hand side yields that $\vec{x}^*(C) \leq \tfrac{1}{2}k =
  \tfrac{1}{2}\card{C}$. Hence, the weight of~$C$ with respect to
  $\tfrac{1}{2} \vones - \vec{x}^*$ can be bounded as follows
  \[
  \sum_{a \in C} \bigl(\tfrac{1}{2} - x_a^*\bigr)  = \tfrac{1}{2} \card{C} - \vec{x}^*(C)
  \geq 0,
  \]
  which proves the lemma.
\end{proof}


We have now reduced the separation problem to finding a shortest cycle
in a weighted undirected graph $G = (V, E)$ without negative cycles.

By using $T$-join techniques, one can compute a shortest path in an
undirected graph without negative cycles in $\order{n_i(m_i + n_i \log
  n_i)}$ time, where in this formula $n_i = \card{\levels^i}$ and $m_i
= \card{A_i}$; see Schrijver~\cite[Chapter 29]{Sch03}. It follows that
a shortest cycle can be computed in $\order{m_i n_i(m_i + n_i \log
  n_i)}$ time. Since $\card{A_i} \leq (i+2) n_i$, this leads to an
$\order{(d+1)^2 n^3 + (d+1) n^3 \log n}$ overall algorithm, where $n
:= \card{\faces}$ is the number of faces and $d$ is the dimension of the
complex.

\subsubsection{Transforming the Graph}
\label{sec:TwoPathGraph}

\begin{figure}
  \psfrag{w1}{\scriptsize $w_1$}
  \psfrag{w2}{\scriptsize $w_2$}
  \psfrag{w3}{\scriptsize $w_3$}
  \psfrag{u1}{\scriptsize $u_1$}
  \psfrag{u2}{\scriptsize $u_2$}
  \psfrag{u3}{\scriptsize $u_3$}
  \psfrag{u4}{\scriptsize $u_4$}
  \includegraphics[width=5cm]{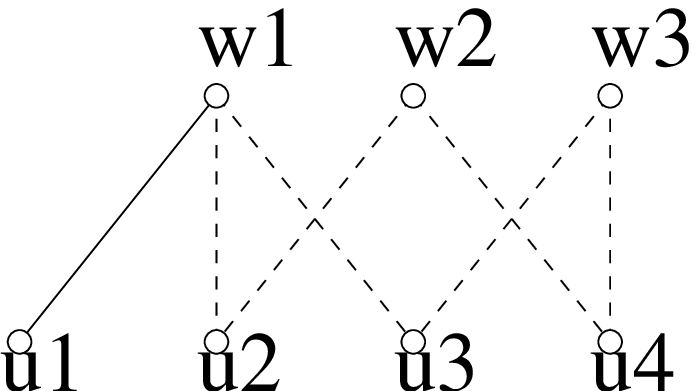}
  \hspace{10ex}
  \psfrag{u1u2w1}[r][r]{\scriptsize $(\{u_1,u_2\},w_1)$}
  \psfrag{u2u3w1}[r][r]{\scriptsize $(\{u_2,u_3\},w_1)$}
  \psfrag{u1u3w1}[r][r]{\scriptsize $(\{u_1,u_3\},w_1)$}
  \psfrag{u2u4w2}{\scriptsize $(\{u_2,u_4\},w_2)$}
  \psfrag{u3u4w3}{\scriptsize $(\{u_3,u_4\},w_3)$}
  \includegraphics[width=7cm]{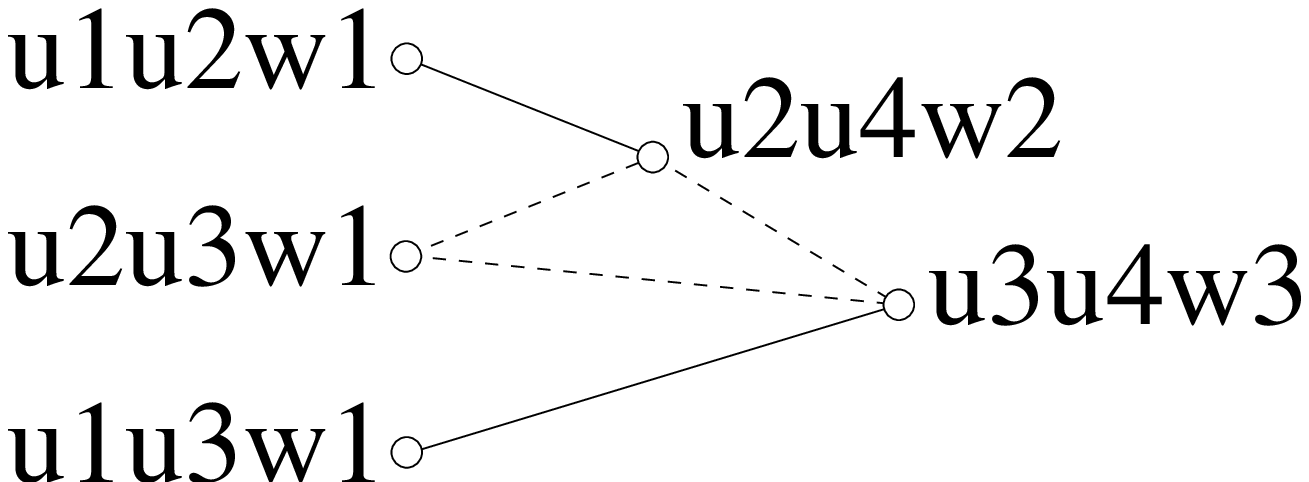}
  \caption{Example of the construction in Section~\ref{sec:TwoPathGraph}. \emph{Left}: original
    graph~$G$. \emph{Right:} constructed graph~$G'$. The $6$-cycle on
    the left corresponds to the $3$-cycle on the right (both shown
    with dashed lines).}
  \label{fig:TwoPath}
\end{figure}

Another method for the separation problem of cycle inequalities, which
is easier to implement, works as follows.

Let $G = (U \,\dot{\cup}\, W,\, E)$ be a bipartite graph, e.g., $G =
H_i$ ($i \in \{0,\dots, d-1\}$), the $i$-th level of the Hasse
diagram. Let $\ell: E \rightarrow \R_{\geq 0}$ be a length function
for the edges of $G$. In the following we write $\ell(u,v) =
\ell(v,u)$ for the length $\ell(\{u,v\})$.

We construct a graph $G' = (V', E')$ and lengths $\ell': E'
\rightarrow \R_{\geq 0}$ as follows; see Figure~\ref{fig:TwoPath} for
an example. The set of nodes of $G'$ is
\[
\big\{ (\{u,u'\}, w) \suchthat \{u, u'\} \subseteq U,\; w \in W,\;
\{u,w\} \in E,\; \{u',w\} \in E \big\}.
\]
We have an edge between two nodes $(\{u_1,u_1'\},w_1)$ and $(\{u_2,
u_2'\}, w_2)$ if
\[
\{u_1,u_1'\} \cap \{u_2,u_2'\} \neq \varnothing \quad \text{and}\quad w_1 \neq w_2.
\]
The length of such an edge $e'$ is defined by
\[
\ell'(e') = \tfrac{1}{2} \big( \ell(u_1,w_1) + \ell(u_1',w_1) +
\ell(u_2,w_2) + \ell(u_2',w_2) \big).
\]
We now consider the relation of cycles in $G$ and $G'$.

\begin{lemma}\label{lemma:Cycles}
  $C = (u_0, w_0, u_1, w_1, \dots, w_{k-1}, u_1)$ is a cycle in
  $G$ with $k > 1$ of length $\ell(C)$ if and only if
  \[
  C' = \big( (\{u_0,u_1\}, w_0), (\{u_1,u_2\},w_1), \dots,
  (\{u_{k-1},u_1\}, w_{k-1}), (\{u_0,u_1\}, w_0) \big)
  \]
  is a cycle in $G'$ with $\ell'(C') = \ell(C)$.
\end{lemma}
\begin{proof}
  In the following we compute with indices modulo $k$.
  
  First observe that $C'$ is well defined: Each $(\{u_i,
  u_{i+1}\},w_i)$ is a node of~$G'$, since $\{u_i, w_i\}$, $\{u_{i+1},
  w_i\} \in E$.  Furthermore,
  \[
  \big\{ (\{u_i,u_{i+1}\}, w_i), (\{u_{i+1},u_{i+2}\},w_{i+1})\big\}
  \in E'
  \]
  since $w_i \neq w_{i+1}$ and $\{u_i,u_{i+1}\} \cap
    \{u_{i+1},u_{i+2}\} \neq \varnothing$ (because $k \geq 1$). It is a
    cycle since $k \geq 2$.

  The weight of $C'$ can be calculated as follows:
  \begin{align*}
    \ell'(C') & = \sum_{i=0}^k \ell\Big( \big\{ (\{u_i,u_{i+1}\}, w_i),
  (\{u_{i+1},u_{i+2}\},w_{i+1})\big\}\Big)\\
  & = \sum_{i=0}^k \tfrac{1}{2} \big( \ell(u_i,w_i) + \ell(u_{i+1},w_i) +
  \ell(u_{i+1},w_{i+1}) + \ell(u_{i+2},w_{i+1}) \big)\\
  & = \sum_{i=0}^k \big( \ell(u_i,w_i) + \ell(u_{i+1},w_i) \big)
  = \ell(C),
  \end{align*}
  where the second to last equation follows since every edge in $C$ occurs
  twice in the summation.
\end{proof}

The previous lemma does not cover cycles in~$G$ of length four.  These
do not occur for the case of $G=H_i$, since $H_i$ is a level in the
Hasse diagram of a \emph{simplicial} complex.  Moreover, cycles of
length four can readily be detected in the construction of~$G'$ and
handled accordingly (there is only a polynomial number of them).

To solve our separation problem, let $G = H_i$, $i \in \{0,\dots,
d-1\}$, and $\ell(e) = x^*_e$ for $e \in G$.  Then we have $\ell'(e')
\in [0,1]$ for each $e' \in E'$, because of the matching inequalities.
We now set $\tilde{\ell}(e') = 1 - \ell'(e')$ for $e' \in G'$ and
hence $\tilde{\ell}(e') \in [0,1]$. Let $C$ be a cycle in $G$ with at
least six edges and $C'$ be the corresponding cycle in $G'$.  Note
that $\card{C'} = \tfrac{1}{2} \card{C}$. We have:
\begin{align*}
  & \tilde{\ell}(C') = \sum_{e' \in C'} \tilde{\ell}(e') = \sum_{e' \in C'} (1 - \ell'(e')) < 1 \\
  & \Leftrightarrow \quad
  \sum_{e' \in C'} \ell'(e') > \card{C'} - 1\\
  & \Leftrightarrow \quad
  \ell'(C') > \card{C'} - 1\\
  & \Leftrightarrow \quad
  \ell(C) > \tfrac{1}{2} \card{C} - 1\qquad\qquad\text{(by Lemma~\ref{lemma:Cycles}).}
\end{align*}
Hence, $C$ violates the cycle inequality~\eqref{eq:cycle_inequalities}
if and only if $\tilde{\ell}(C') < 1$. Since $\tilde{\ell}(e') \geq
0$, we can use the Floyd-Warshall algorithm to solve the separation
problem in time $\order{\card{V'}^3}$; see Korte and
Vygen~\cite{KorV02}.

If $G = H_i$ and $W$ is the part arising from the higher dimensional
faces, we have $\card{V'} = \tbinom{i+2}{2} \card{W} = \tbinom{i+2}{2}
f_{i+1}$. This leads to an $\order{(d+1)^6 n^3}$ algorithm for
separating cycle inequalities, which is roughly as fast as the method
discussed in Section~\ref{sec:undirected_shortest_paths}, but much
easier to implement.

\section{Computational Results}

In this section we report on computational experience with a
branch-and-cut algorithm along the lines of
Section~\ref{sec:ip_formulation}. The C++ implementation uses the
framework SCIP (Solving Constraint Integer Programs) by Achterberg,
see~\cite{Ach04}. It furthermore builds on \texttt{polymake}; see
\cite{GawJ00, GawJ04}.  As an LP solver we used CPLEX 9.0.

As the basis of our implementation we take the formulation of \MMM in
Section~\ref{sec:ip_formulation}. Matching
inequalities~\eqref{eq:matching_inequalities} and Betti
inequalities~\eqref{eq:betti_inequality} (together with variable
bounds) form the initial LP. Cycle
inequalities~\eqref{eq:cycle_inequalities} are separated as described
in Section~\ref{sec:TwoPathGraph}. Additionally, Gomory cuts are
added. As a branching rule we use \emph{reliability branching}
implemented in SCIP, a variable branching rule introduced by
Achterberg, Koch, and Martin~\cite{AchKM04}.


We implemented the following primal heuristic. First a simple greedy
algorithm is run: We start with the empty matching $M = \varnothing$.
We add arcs of the Hasse diagram to $M$ in the order of decreasing
value of the current LP solution as long as $M$ stays an acyclic
matching (which can easily be tested). Then the outcome is iteratively
improved by a method described in Forman~\cite{For98}: One searches
for a unique path between two critical faces in~$\HD(M)$. Such a path
is alternating with respect to~$M$. Then~$M$ can be augmented along
the path (the new matching is the symmetric difference of~$M$ and the
path). As is easily seen, this generates an acyclic matching, because
the path is unique. This heuristic turns out to be extremely
successful; see below.

\begin{table}[tb]
  \caption{Computational results of the branch-and-cut
    algorithm with separating cycle inequalities and Gomory cuts.}\label{tab:results1}
  \begin{center}
    \renewcommand{\arraystretch}{0.9}
    \begin{tabular*}{\linewidth}{@{\extracolsep{\fill}}lrrrrrrrr@{}}\toprule
      name            &    n &     m & d &  nodes &depth & time    & $\beta$ &  c  \\\midrule
      solid\_2\_torus &   24 &    42 & 2 &      1 &    0 &    0.00 &  2 &  2  \\
      simon2          &   31 &    60 & 2 &      1 &    0 &    0.00 &  1 &  1  \\
      projective      &   31 &    60 & 2 &      1 &    0 &    0.01 &  3 &  3  \\
      bjorner         &   32 &    63 & 2 &      1 &    0 &    0.05 &  2 &  2  \\
      nonextend       &   39 &    77 & 2 &      6 &    5 &    0.16 &  1 &  1  \\
      simon           &   41 &    82 & 2 &      1 &    0 &    0.18 &  1 &  1  \\
      dunce           &   49 &    99 & 2 &    385 &   10 &    2.62 &  1 &  3  \\
      c-ns3           &   63 &   128 & 2 &    349 &   10 &    3.47 &  1 &  3  \\
      c-ns            &   75 &   152 & 2 &     28 &   10 &    1.95 &  1 &  3  \\
      c-ns2           &   79 &   159 & 2 &     14 &    7 &    1.11 &  1 &  1  \\
      ziegler         &  119 &   310 & 3 &      1 &    0 &    0.01 &  1 &  1  \\
      gruenbaum       &  167 &   434 & 3 &      1 &    0 &   25.24 &  1 &  1  \\
      lockeberg       &  216 &   600 & 3 &      1 &    0 &   36.25 &  2 &  2  \\
      rudin           &  215 &   578 & 3 &     77 &   30 &  103.78 &  1 &  1  \\
      mani-walkup-D   &  392 &  1112 & 3 &    111 &   23 &  512.81 &  2 &  2 \\
      mani-walkup-C   &  464 &  1312 & 3 &    135 &   83 & 1658.02 &  2 &  2 \\\addlinespace

      MNSB            &  103 &   267 & 3 &     12 &   10 &   73.39 &  1 &  1 \\
      MNSS            &  250 &   698 & 3 &    292 &  110 &  750.36 &  2 &  2 \\
      CP2             &  255 &   864 & 4 &    230 &   80 &  558.14 &  3 &  3 \\\bottomrule
    \end{tabular*}
  \end{center}
\end{table}

We tested the implementation on a set of simplicial complexes
collected by Hachimori; see~\cite{Hac01} for more details.
Additionally, we considered the following complexes: \texttt{CP2}
(complex projective plane), \texttt{CP2+CP2} (connected sum of
\texttt{CP2} with itself), \texttt{MNSB} and \texttt{MNSS} ((vertex)
minimal non-shellable ball and sphere, respectively; see
Lutz~\cite{Lut04}).

All computational experiments were run on a 3 GHz Pentium machine
running Linux. In the tables of computational results, $n$ denotes the
number of faces, $m$ the number of arcs in the Hasse diagram (= number
of variables), $d$ the dimension, \emph{nodes} the number of nodes in
the branch-and-bound tree, \emph{depth} the maximal depth in the tree,
\emph{time} the computation time in seconds, $\beta$ is the lower
bound obtained by adding all Betti
inequalities~\eqref{eq:betti_inequality}, and $c$ the number of
critical faces in the optimal solution.

Our implementation could not solve the larger problems of Hachimori's
collection in reasonable time: \texttt{bing}, \texttt{knot},
\texttt{poincare}, \texttt{nonpl\_sphere}, and \texttt{nc\_sphere}. In
fact, for \texttt{poincare} we ran our code in different settings,
each for about a week -- without success.

Table~\ref{tab:results1} shows the results of a computation where we
separate cycle inequalities and Gomory cuts and run the heuristic
every 10th level. At most seven separation rounds of cycle
inequalities were performed at a node. We do not report results on the
problems by Moriyama and Takeuchi in Hachimori's collection -- they
all could be solved within a second. The version with cut separation
could not solve \texttt{CP2+CP2} within 90 minutes.

\begin{table}[tb]
  \caption{Computational results of the branch-and-cut
    algorithm without separation.}\label{tab:results2}
  \begin{center}
    \renewcommand{\arraystretch}{0.9}
    \begin{tabular*}{\linewidth}{@{\extracolsep{\fill}}lrrrrrrrr@{}}\toprule
      name            &    n &     m & d &  nodes &depth & time    & $\beta$ &  c  \\\midrule
      solid\_2\_torus &   24 &    42 & 2 &      1 &    0 &    0.00 &  2 &  2  \\
      simon2          &   31 &    60 & 2 &      1 &    0 &    0.01 &  1 &  1  \\
      projective      &   31 &    60 & 2 &      1 &    0 &    0.00 &  3 &  3  \\
      bjorner         &   32 &    63 & 2 &      1 &    0 &    0.01 &  2 &  2  \\
      nonextend       &   39 &    77 & 2 &      3 &    2 &    0.02 &  1 &  1  \\
      simon           &   41 &    82 & 2 &      4 &    3 &    0.02 &  1 &  1  \\
      dunce           &   49 &    99 & 2 & 168367 &   42 &  145.60 &  1 &  3  \\
      c-ns3           &   63 &   128 & 2 & 3665581 &  53 & 3940.40 &  1 &  3  \\
      c-ns            &   75 &   152 & 2 & 16625713 & 58 &19359.69 &  1 &  3  \\
      c-ns2           &   79 &   159 & 2 &      4 &    3 &    0.03 &  1 &  1  \\
      ziegler         &  119 &   310 & 3 &      1 &    0 &    0.01 &  1 &  1  \\
      gruenbaum       &  167 &   434 & 3 &     21 &   20 &    0.68 &  1 &  1  \\
      lockeberg       &  216 &   600 & 3 &      1 &    0 &    0.05 &  2 &  2  \\
      rudin           &  215 &   578 & 3 &     81 &   80 &    3.18 &  1 &  1  \\
      mani-walkup-D   &  392 &  1112 & 3 &    107 &  100 &    2.00 &  2 &  2  \\
      mani-walkup-C   &  464 &  1312 & 3 &   1498 &  456 &   30.54 &  2 &  2  \\\addlinespace
      MNSB            &  103 &   267 & 3 &      1 &    0 &    0.01 &  1 &  1  \\
      MNSS            &  250 &   698 & 3 &    163 &  126 &    4.63 &  2 &  2  \\
      CP2             &  255 &   864 & 4 &    198 &  190 &    4.77 &  3 &  3  \\
      CP2+CP2         &  460 &  1592 & 4 &   5178 &  534 &  110.21 &  4 &  4  \\\bottomrule
    \end{tabular*}
  \end{center}
\end{table}

For most problems the bound obtained by adding Betti
inequalities~\eqref{eq:betti_inequality}, as indicated in
column~``$\beta$'', is tight. This means that the algorithm is done
once an optimal solution is found. This usually happens very fast and
shows that the heuristic is efficient. In fact, there are only three
problems for which the bound is not tight and could be solved by our
algorithm (\texttt{dunce}, \texttt{c-ns}, and \texttt{c-ns3}). These
three problems are solved easily by the version with cut separation.
In our problem set there exists no hard but still solvable problem
with a ``Betti bound'' which is not sharp. We can therefore not
estimate the limits of our implementation for these cases
(\texttt{poincare} is the next larger problem of this kind with 1112
variables, but we could not solve it).

The tractability of problems with a tight ``Betti bound'' is supported
by the results obtained by running the implementation without any
separation; see Table~\ref{tab:results2}.  Only integer solutions are
checked whether they are acyclic and the heuristic is run every 10th
level. This essentially is a test of the performance of the primal
heuristic. Indeed, all problems with tight ``Betti bound'' were solved
within a few seconds (\texttt{CP2+CP2} and \texttt{mani-walkup-C}
being the exception, but could be solved within two minutes). The
results for the problems \texttt{c-ns}, \texttt{c-ns3}, and
\texttt{dunce} show that the cycle inequalities and Gomory cuts are
very effective in reducing the number of nodes in the tree and the
computing time for problems where the ``Betti bound'' is not sharp.

Summarizing, we can say that our implementation can solve large
instances with up to about 1500 variables if the bounds from the Betti
numbers are tight and small instances with up to about 150 variables if
the bounds are not tight.  In all the instances computed so far, the
topology of the spaces involved was known.  In the future, we plan to
apply our techniques to other cases.

\subsection*{Acknowledgment}

We are indebted to Tobias Achterberg for his support of the implementation.

\bibliographystyle{siam}
\bibliography{morse}

\end{document}